\documentclass[11pt,twoside]{article}
\usepackage{amsmath, amsthm, amscd, amsfonts, amssymb, graphicx, color}

\date{}

\begin{document}

\centerline{}

\centerline {\Large{\bf  Constructions of bi-$g$-fusion frame in Hilbert space}}

\centerline{}

%% My definition
\newcommand{\mvec}[1]{\mbox{\bfseries\itshape #1}}
\centerline{}
\centerline{\textbf{Prasenjit Ghosh}}
\centerline{Department of Mathematics,  }
\centerline{Barwan N.S. High School(HS),}
\centerline{Barwan, Murshidabad, 742132, West Bengal, India}
\centerline{e-mail: prasenjitpuremath@gmail.com}
\centerline{}
\centerline{\textbf{T. K. Samanta}}
\centerline{Department of Mathematics, Uluberia College,}
\centerline{Uluberia, Howrah, 711315,  West Bengal, India}
\centerline{e-mail: mumpu$_{-}$tapas5@yahoo.co.in}

\newtheorem{Theorem}{\quad Theorem}[section]

\newtheorem{definition}[Theorem]{\quad Definition}

\newtheorem{theorem}[Theorem]{\quad Theorem}

\newtheorem{remark}[Theorem]{\quad Remark}

\newtheorem{corollary}[Theorem]{\quad Corollary}

\newtheorem{note}[Theorem]{\quad Note}

\newtheorem{lemma}[Theorem]{\quad Lemma}

\newtheorem{example}[Theorem]{\quad Example}

\newtheorem{result}[Theorem]{\quad Result}
\newtheorem{conclusion}[Theorem]{\quad Conclusion}

\newtheorem{proposition}[Theorem]{\quad Proposition}

\begin{abstract}
\textbf{\emph{The concept of a bi-$g$-fusion frame for a Hilbert space, which is a generalizations of a controlled $g$-fusion frame, is introduced and an example is given.\,Finally,  bi-$g$-fusion frame in tensor product of Hilbert spaces is considered.}}
\end{abstract}
{\bf Keywords:}  \emph{Frame, $g$-fusion frame, bi-frame, Tensor product of Hilbert spaces.}

{\bf 2010 Mathematics Subject Classification:} \emph{42C15; 46C07; 46M05; 47A80.}

%=====================================
\section{Introduction}
%=====================================

In 1952, Duffin and Schaeffer \cite{Duffin} introduced Frame for Hilbert space to study some fundamental problems in non-harmonic Fourier series.\;The formal definition of frame in the abstract Hilbert spaces was given by Daubechies et al.\,\cite{Daubechies} in 1986.\;Frame theory has been widely used in signal and image processing, filter bank theory, coding and communications, system modeling and so on.\;Several generalizations of frames  namely, \,$g$-frames \cite{Sun}, fusion frames \cite{Kutyniok}, \,$g$-fusion frames \cite{Ahmadi} etc.\;had been introduced in recent times.\,Generalized atomic subspaces for operators in Hilbert spaces were studied by P.\,Ghosh and T.\,K.\,Samanta {\cite{Ghosh}} and they were also presented the stability of dual \,$g$-fusion frames in Hilbert spaces in {\cite{P}}.\,Biframe is also a generalization of controlled frame in Hilbert space which was studied by M. F. Parizi et al.\,\cite{MF}.\,To define frame in Hilbert space, only one sequence is needed, but for a biframe, two sequences are needed.\,A pair of sequences \,$\left(\,\left\{\,f_{i}\,\right\}_{i \,=\, 1}^{\,\infty}\,,\, \left\{\,g_{i}\,\right\}_{i \,=\, 1}^{\,\infty}\,\right)$\, in \,$H$\, is called a biframe for \,$H$\, if there exist positive constants \,$A$\, and \,$B$\, such that
\[A\; \|\,f\,\|^{\,2} \,\leq\, \sum\limits_{i \,=\, 1}^{\infty}\, \left <\,f\,,\, f_{\,i} \, \right >\,\left<\,g_{\,i}\,,\, f\,\right> \,\leq\, B \,\|\,f\,\|^{\,2}\; \;\forall\; f \,\in\, H.\] 
The constants \,$A$\, and \,$B$\, are called lower and upper biframe bounds, respectively.

The basic concepts of tensor product of Hilbert spaces were presented by S.\,Rabin-son \cite{S}.\;Frames and Bases in Tensor Product of Hilbert spaces were introduced by A.\,Khosravi and M.\,S.\,Asgari \cite{A}.\;Reddy et al.\,\cite{Upender} also studied the frame in tensor product of Hilbert spaces and presented the tensor frame operator on tensor product of Hilbert spaces.

In this paper, the idea of a bi $g$-fusion frame in Hilbert space is presented and some of their properties are going to be established.\,We also given an example of a bi $g$-fusion frame in Hilbert space.\,Bi $g$-fusion frame in tensor product of Hilbert spaces is discussed and the relation between the frame operators for the pair of \,$g$-fusion Bessel sequences in Hilbert spaces and their tensor product are obtained.

Throughout this paper,\;$H$\, is considered to be separable Hilbert space with associated inner products \,$\left <\,\cdot \,,\, \cdot\,\right>$\,.\,$I_{H}$\, denotes the identity operator on \,$H$.\,$\mathcal{B}\,(\,H\,)$\; denote the space of all bounded linear operators on \,$H$.\,$P_{\,V}$\, denote the orthogonal projection onto the closed subspace \,$V \,\subset\, H$.\;$\left\{\,V_{i}\,\right\}_{ i \,\in\, I}$\, is the collections of closed subspaces of \,$H$, where \,$I$\, is subset of  integers \,$\mathbb{Z}$.\;$\left\{\,H_{i}\,\right\}_{ i \,\in\, I}$\, is the collection of Hilbert spaces.\,$\left\{\,\Lambda_{i} \,\in\, \mathcal{B}\,(\,H \,,\, H_{i}\,)\,\right\}_{i \,\in\, I}$\, denotes the sequence of operators.\;Define the space
\[l^{\,2}\left(\,\left\{\,H_{i}\,\right\}_{ i \,\in\, I}\,\right) \,=\, \left \{\,\{\,f_{\,i}\,\}_{i \,\in\, I} \,:\, f_{\,i} \;\in\; H_{i},\; \sum\limits_{\,i \,\in\, I}\, \left \|\,f_{\,i}\,\right \|^{\,2} \,<\, \infty \,\right\}\]
with inner product is given by \,$\left<\,\{\,f_{\,i}\,\}_{i \,\in\, I} \,,\, \{\,g_{\,i}\,\}_{ i \,\in\, I}\,\right> \;=\; \sum\limits_{\,i \,\in\, I}\, \left<\,f_{\,i} \,,\, g_{\,i}\,\right>_{H_{i}}$.\;Clearly \,$l^{\,2}\left(\,\left\{\,H_{i}\,\right\}_{ i \,\in\, I}\,\right)$\; is a Hilbert space with respect to the above inner product \cite{Ahmadi}.

%=====================================
\section{Preliminaries}
%=====================================

\smallskip\hspace{.6 cm}

\begin{theorem}\cite{Gavruta}\label{th0.001}
Let \,$V \,\subset\, H$\; be a closed subspace and \,$T \,\in\, \mathcal{B}\,(\,H\,)$.\;Then \,$P_{\,V}\, T^{\,\ast} \,=\, P_{\,V}\,T^{\,\ast}\, P_{\,\overline{T\,V}}$.\;If \,$T$\; is an unitary operator (\,i\,.\,e., \,$T^{\,\ast}\, T \,=\, I_{H}$\,), then \,$P_{\,\overline{T\,V}}\;T \,=\, T\,P_{\,V}$.
\end{theorem}

\begin{theorem}\cite{Jain}\label{th0.01}
The set \,$\mathcal{S}\,(\,H\,)$\; of all self-adjoint operators on \,$H$\; is a partially ordered set with respect to the partial order \,$\leq$\, which is defined as for \,$T,\,S \,\in\, \mathcal{S}\,(\,H\,)$ 
\[T \,\leq\, S \,\Leftrightarrow\, \left<\,T\,f \,,\, f\,\right> \,\leq\, \left<\,S\,f \,,\, f\,\right>\; \;\forall\; f \,\in\, H.\] 
\end{theorem}

\begin{theorem}(\,Douglas' factorization theorem\,)\,{\cite{Douglas}}\label{th1}
Let \;$U,\, V \,\in\, \mathcal{B}\,(\,H\,)$.\,Then the following conditions are equivalent:
\begin{itemize}
\item[(\,I\,)]\;\;$\mathcal{R}\,(\,U\,) \,\subseteq\, \mathcal{R}\,(\,V\,)$.
\item[(\,II\,)]\;\;$U\, U^{\,\ast} \,\leq\, \lambda^{\,2}\; V\,V^{\,\ast}$\; for some \,$\lambda \,>\, 0$.
\item[(\,III\,)]\;\;$U \,=\, V\,W$\, for some bounded linear operator \,$W$\, on \,$H$.
\end{itemize}
\end{theorem}

\begin{definition}\cite{Ahmadi}
Let \,$\left\{\,v_{i}\,\right\}_{ i \,\in\, I}$\, be a collection of positive weights.\;Then the family \,$\Lambda \,=\, \{\,\left(\,V_{i},\, \Lambda_{i},\, v_{i}\,\right)\,\}_{i \,\in\, I}$\; is called a generalized fusion frame or a g-fusion frame for \,$H$\; respect to \,$\left\{\,H_{i}\,\right\}_{i \,\in\, I}$\; if there exist constants \,$0 \,<\, A \,\leq\, B \,<\, \infty$\; such that
\begin{equation}\label{eq1.01}
A \;\left \|\,f \,\right \|^{\,2} \,\leq\, \sum\limits_{\,i \,\in\, I}\,v_{i}^{\,2}\, \left\|\,\Lambda_{i}\,P_{\,V_{i}}\,(\,f\,) \,\right\|^{\,2} \,\leq\, B \; \left\|\, f \, \right\|^{\,2}\; \;\forall\; f \,\in\, H.
\end{equation}
The constants \,$A$\; and \,$B$\; are called the lower and upper bounds of g-fusion frame, respectively.\,If \,$A \,=\, B$\; then \,$\Lambda$\; is called tight g-fusion frame and if \;$A \,=\, B \,=\, 1$\, then we say \,$\Lambda$\; is a Parseval g-fusion frame.\;If  \,$\Lambda$\; satisfies the inequality
\[\sum\limits_{\,i \,\in\, I}\,v_{i}^{\,2}\, \left\|\,\Lambda_{i}\,P_{\,V_{i}}\,(\,f\,) \,\right\|^{\,2} \,\leq\, B \; \left\|\, f \, \right\|^{\,2}\; \;\forall\; f \,\in\, H\]
then it is called a g-fusion Bessel sequence with bound \,$B$\; in \,$H$. 
\end{definition}

\begin{definition}\cite{Ahmadi}\label{defn1}
Let \,$\Lambda \,=\, \{\,\left(\,V_{i},\, \Lambda_{i},\, v_{i}\,\right)\,\}_{i \,\in\, I}$\, be a g-fusion Bessel sequence in \,$H$\, with a bound \,$B$.\;The synthesis operator \,$T_{\Lambda}$\, of \,$\Lambda$\; is defined as 
\[ T_{\Lambda} \,:\, l^{\,2}\left(\,\left\{\,H_{i}\,\right\}_{ i \,\in\, I}\,\right) \,\to\, H,\]
\[T_{\Lambda}\,\left(\,\left\{\,f_{\,i}\,\right\}_{i \,\in\, I}\,\right) \,=\,  \sum\limits_{\,i \,\in\,I}\, v_{i}\, P_{\,V_{i}}\,\Lambda_{i}^{\,\ast}\,f_{i}\; \;\;\forall\; \{\,f_{i}\,\}_{i \,\in\, I} \,\in\, l^{\,2}\left(\,\left\{\,H_{i}\,\right\}_{ i \,\in\, I}\,\right)\] and the analysis operator is given by 
\[ T_{\Lambda}^{\,\ast} \,:\, H \,\to\, l^{\,2}\left(\,\left\{\,H_{i}\,\right\}_{ i \,\in\, I}\,\right),\; T_{\Lambda}^{\,\ast}\,(\,f\,) \,=\,  \left\{\,v_{i}\,\Lambda_{i}\, P_{\,V_{i}}\,(\,f\,)\,\right\}_{ i \,\in\, I}\; \;\forall\; f \,\in\, H.\]
The g-fusion frame operator \,$S_{\Lambda} \,:\, H \,\to\, H$\; is defined as follows:
\[S_{\Lambda}\,(\,f\,) \,=\, T_{\Lambda}\,T_{\Lambda}^{\,\ast}\,(\,f\,) \,=\, \sum\limits_{\,i \,\in\, I}\, v_{i}^{\,2}\; P_{\,V_{i}}\, \Lambda_{i}^{\,\ast}\; \Lambda_{i}\, P_{\,V_{i}}\,(\,f\,)\; \;\forall\; f \,\in\, H.\]
\end{definition}

\begin{note}\cite{Ahmadi}
Let \,$\Lambda \,=\, \{\,\left(\,V_{i},\, \Lambda_{i},\, v_{i}\,\right)\,\}_{i \,\in\, I}$\, be a g-fusion Bessel sequence in \,$H$.\,Then it can be easily verify that for all \,$f \,\in\, H$ 
\[\left<\,S_{\Lambda}\,(\,f\,) \,,\, f\,\right> \,=\, \sum\limits_{\,i \,\in\, I}\, v_{i}^{\,2}\, \left\|\,\Lambda_{i}\, P_{\,V_{i}}\,(\,f\,) \,\right\|^{\,2} \,=\, \left\|\,T_{\Lambda}^{\,\ast}\,(\,f\,)\,\right\|^{\,2}.\]
If \,$\Lambda$\, is a g-fusion frame with bounds \,$A$\, and \,$B$\, then from (\ref{eq1.01}),
\[\left<\,A\,f \,,\, f\,\right> \,\leq\, \left<\,S_{\Lambda}\,(\,f\,) \,,\, f\,\right> \,\leq\, \left<\,B\,f \,,\, f\,\right>\; \;\forall\; f \,\in\, H.\]
Now, according to the Theorem (\ref{th0.01}), we can write, \,$A\,I_{H} \,\leq\,S_{\Lambda} \,\leq\, B\,I_{H}$.\;The operator \,$S_{\Lambda}$\; is bounded, self-adjoint, positive and invertible.\;Also, \,$ B^{\,-1}\,I_{H} \,\leq\, S_{\,\Lambda}^{\,-1} \,\leq\, A^{\,-1}\,I_{H}$.\;Hence, reconstruction formula for any \,$f \,\in\, H$, is given by  
\[ f \,=\, \sum\limits_{\,i \,\in\, I}\,v^{\,2}_{\,i}\,P_{\,V_{i}}\,\Lambda^{\,\ast}_{i}\,\Lambda_{i}\,P_{\,V_{i}}\,S^{\,-\, 1}_{\Lambda}\,(\,f\,) \,=\, \sum\limits_{\,i \,\in\, I}\,v^{\,2}_{\,i}\,S^{\,-\, 1}_{\Lambda}\,P_{\,V_{i}}\,\Lambda^{\,\ast}_{i}\,\Lambda_{i}\,P_{\,V_{i}}\,(\,f\,).\]
\end{note}

\begin{definition}\cite{Ahmadi}\label{defn2}
Let \,$\Lambda \,=\, \{\,\left(\,V_{i},\, \Lambda_{i},\, v_{i}\,\right)\,\}_{i \,\in\, I}$\, be a g-fusion frame for \,$H$\, with frame operator \,$S_{\Lambda}$.\;Then the g-fusion frame \,$\left\{\,\left(\,S^{\,-\, 1}_{\Lambda}\,V_{i},\, \Lambda_{i}\,P_{\,V_{i}}\,S^{\,-\, 1}_{\Lambda},\, v_{i}\,\right)\,\right\}_{i \,\in\, I}$\, is called the canonical dual g-fusion frame of \,$\Lambda$. 
\end{definition}

\begin{definition}\cite{Ghosh}
Let \,$\Lambda \,=\, \left\{\,\left(\,V_{i},\, \Lambda_{i},\, v_{i}\,\right)\,\right\}_{i \,\in\, I}$\, and \,$\Lambda^{\,\prime} \,=\, \left\{\,\left(\,V^{\,\prime}_{i},\, \Lambda^{\,\prime}_{i},\, v^{\,\prime}_{i}\,\right)\,\right\}_{i \,\in\, I}$\, be two g-fusion Bessel sequences in \,$H$\; with bounds \,$D_{\,1}$\, and \,$D_{\,2}$, respectively.\;Then the operator \,$S_{\Lambda\,\Lambda^{\,\prime}} \,:\, H \,\to\, H$, defined by
\[S_{\Lambda\,\Lambda^{\,\prime}}\,(\,f\,) \;=\; \sum\limits_{\,i \,\in\, I}\,v_{\,i}\,v^{\,\prime}_{\,i}\,P_{\,V_{i}}\,\Lambda_{i}^{\,\ast}\,\Lambda_{i}^{\,\prime} \,P_{\,V^{\,\prime}_{i}}\,(\,f\,)\;\; \;\forall\; f \,\in\, H,\] is called the frame operator for the pair of g-fusion Bessel sequences \,$\Lambda$\, and \,$\Lambda^{\,\prime}$.
\end{definition}

\begin{definition}\cite{HS}
Let \,$\left\{\,W_{j}\,\right\}_{ j \,\in\, J}$\, be a collection of closed subspaces of \,$H$\, and \,$\left\{\,v_{j}\,\right\}_{ j \,\in\, J}$\, be a collection of positive weights.\,Let \,$\left\{\,H_{j}\,\right\}_{ j \,\in\, J}$\, be a sequence of Hilbert spaces, \,$T,\, U \,\in\, \mathcal{G}\,\mathcal{B}\,(\,H\,)$\, and \,$\Lambda_{j} \,\in\, \mathcal{B}\,(\,H,\, H_{j}\,)$\, for each \,$j \,\in\, J$.\,Then the family \,$\Lambda_{T\,U} \,=\, \left\{\,\left(\,W_{j},\, \Lambda_{j},\, v_{j}\,\right)\,\right\}_{j \,\in\, J}$\, is a \,$(\,T,\,U\,)$-controlled $g$-fusion frame for \,$H$\, if there exist constants \,$0 \,<\, A \,\leq\, B \,<\, \infty$\, such that 
\begin{equation}\label{eqn1.1}
A\,\|\,f\,\|^{\,2} \,\leq\, \sum\limits_{\,j \,\in\, J}\, v^{\,2}_{j}\,\left<\,\Lambda_{j}\,P_{\,W_{j}}\,U\,f,\,  \Lambda_{j}\,P_{\,W_{j}}\,T\,f\,\right> \,\leq\, \,B\,\|\,f \,\|^{\,2}\; \;\forall\; f \,\in\, H.
\end{equation}
If \,$A \,=\, B$\, then \,$\Lambda_{T\,U}$\, is called \,$(\,T,\,U\,)$-controlled tight g-fusion frame and if \,$A \,=\, B \,=\, 1$\, then we say \,$\Lambda_{T\,U}$\, is a \,$(\,T,\,U\,)$-controlled Parseval g-fusion frame.\,If \,$\Lambda_{T\,U}$\, satisfies only the right inequality of (\ref{eqn1.1}) it is called a \,$(\,T,\,U\,)$-controlled g-fusion Bessel sequence in \,$H$.    
\end{definition}

There are several ways to introduced the tensor product of Hilbert spaces.\;The tensor product of Hilbert spaces \,$H$\, and \,$K$\, is a certain linear space of operators which was represented by Folland in \cite{Folland}, Kadison and Ringrose in \cite{Kadison}.\\

\begin{definition}\cite{Upender}
The tensor product of Hilbert spaces \,$H$\, and \,$K$\, is denoted by \,$H \,\otimes\, K$\, and it is defined to be an inner product space associated with the inner product
\begin{equation}\label{eq1.001}   
\left<\,f \,\otimes\, g \,,\, f^{\,\prime} \,\otimes\, g^{\,\prime}\,\right> \,=\, \left<\,f  \,,\, f^{\,\prime}\,\right>_{\,1}\;\left<\,g  \,,\, g^{\,\prime}\,\right>_{\,2}\; \;\forall\; f,\, f^{\,\prime} \,\in\, H\; \;\text{and}\; \;g,\, g^{\,\prime} \,\in\, K.
\end{equation}
The norm on \,$H \,\otimes\, K$\, is given by 
\begin{equation}\label{eq1.0001}
\left\|\,f \,\otimes\, g\,\right\| \,=\, \|\,f\,\|_{\,1}\;\|\,g\,\|_{\,2}\; \;\forall\; f \,\in\, H\; \;\&\; \,g \,\in\, K.
\end{equation}
The space \,$H \,\otimes\, K$\, is complete with respect to the above inner product.\;Therefore the space \,$H \,\otimes\, K$\, is a Hilbert space.     
\end{definition} 

For \,$Q \,\in\, \mathcal{B}\,(\,H\,)$\, and \,$T \,\in\, \mathcal{B}\,(\,K\,)$, the tensor product of operators \,$Q$\, and \,$T$\, is denoted by \,$Q \,\otimes\, T$\, and defined as 
\[\left(\,Q \,\otimes\, T\,\right)\,A \,=\, Q\,A\,T^{\,\ast}\; \;\forall\; \;A \,\in\, H \,\otimes\, K.\]
It can be easily verified that \,$Q \,\otimes\, T \,\in\, \mathcal{B}\,(\,H \,\otimes\, K\,)$\, \cite{Folland}.\\

\begin{theorem}\cite{Folland}\label{th1.1}
Suppose \,$Q,\, Q^{\prime} \,\in\, \mathcal{B}\,(\,H\,)$\, and \,$T,\, T^{\prime} \,\in\, \mathcal{B}\,(\,K\,)$, then \begin{itemize}
\item[(I)]\hspace{.2cm} \,$Q \,\otimes\, T \,\in\, \mathcal{B}\,(\,H \,\otimes\, K\,)$\, and \,$\left\|\,Q \,\otimes\, T\,\right\| \,=\, \|\,Q\,\|\; \|\,T\,\|$.
\item[(II)]\hspace{.2cm} \,$\left(\,Q \,\otimes\, T\,\right)\,(\,f \,\otimes\, g\,) \,=\, Q\,(\,f\,) \,\otimes\, T\,(\,g\,)$\, for all \,$f \,\in\, H,\, g \,\in\, K$.
\item[(III)]\hspace{.2cm} $\left(\,Q \,\otimes\, T\,\right)\,\left(\,Q^{\,\prime} \,\otimes\, T^{\,\prime}\,\right) \,=\, (\,Q\,Q^{\,\prime}\,) \,\otimes\, (\,T\,T^{\,\prime}\,)$. 
\item[(IV)]\hspace{.2cm} \,$Q \,\otimes\, T$\, is invertible if and only if \,$Q$\, and \,$T$\, are invertible, in which case \,$\left(\,Q \,\otimes\, T\,\right)^{\,-\, 1} \,=\, \left(\,Q^{\,-\, 1} \,\otimes\, T^{\,-\, 1}\,\right)$.
\item[(V)]\hspace{.2cm} \,$\left(\,Q \,\otimes\, T\,\right)^{\,\ast} \,=\, \left(\,Q^{\,\ast} \,\otimes\, T^{\,\ast}\,\right)$.    
\end{itemize}
\end{theorem}

%=====================================
\section{Bi-$g$-fusion frame in Hilbert space}
%=====================================

\smallskip\hspace{.6 cm} In this section, we discuss bi-$g$-fusion frame and the canonical dual bi-$g$-fusion frame in a Hilbert space. 

\begin{definition}
Let \,$V \,=\, \left\{\,V_{\,i}\,\right\}_{\, i \,\in\, I},\; \;W \,=\, \left\{\,W_{\,i}\,\right\}_{ i \,\in\, I}$\, be two families of closed subspaces of \,$H$, \,$\left\{\,v_{\,i}\,\right\}_{\, i \,\in\, I}$\, be family of positive weights i\,.\,e., \,$v_{\,i} \,>\, 0\, \;\forall\; i \,\in\, I$\, and \,$\Lambda_{i}\,,\, \Gamma_{i} \,\in\, \mathcal{B}\,(\,H,\, H_{i}\,)$, for each \,$i \,\in\, I$.\,Then the pair given by \,$\left(\,\Lambda,\, \Gamma\,\right) \,=\, \left(\,\left\{\,\left(\,V_{i},\, \Lambda_{i},\, v_{i}\,\right)\,\right\}_{i \,\in\, I}\,,\,\,\left\{\,\left(\,W_{i},\, \Gamma_{i},\, v_{i}\,\right)\,\right\}_{i \,\in\, I}\,\right)$\, is said to be a bi-$g$-fusion frame for \,$H$\, if there exist constants \,$0 \,<\, A \,\leq\, B \,<\, \infty$\, such that 
\begin{equation}\label{eqn1.1}
A\,\|\,f\,\|^{\,2} \,\leq\, \sum\limits_{\,i \,\in\, I}\, v^{\,2}_{i}\,\left<\,\Lambda_{i}\,P_{\,V_{i}}\,f,\, \Gamma_{i}\,P_{\,W_{i}}\,f\,\right> \,\leq\, \,B\,\|\,f \,\|^{\,2},
\end{equation}
for all \,$f \,\in\, H$,\, where \,$P_{\,V_{i}}$\, and \,$P_{\,W_{i}}$\, are the orthogonal projections of \,$H$\, onto \,$V_{\,i}$\, and \,$W_{\,i}$, respectively.\;The constants \,$A$\, and \,$B$\, are called the frame bounds of \,$\left(\,\Lambda,\,  \Gamma\,\right)$.\;If \,$A \,=\, B$\, then it is called a tight bi-$g$-fusion frame.\;If the pair \,$\left(\,\Lambda,\,  \Gamma\,\right)$\, satisfies only the right inequality of (\ref{eqn1.1}), then it is called a bi-$g$-fusion Bessel sequence in \,$H$\, with bound \,$B$.   
\end{definition}

\begin{remark}
\begin{itemize}
\item[$(I)$]In particular, if \,$\Lambda_{i} \,=\, \Gamma_{i}$\, and \,$V_{i} \,=\, W_{i}$, for each \,$i \,\in\, I$, i.\,e., if \,$\left(\,\Lambda,\,  \Lambda\,\right)$\, is a bi-$g$-fusion frame for \,$H$, then \,$\Lambda$\, is a \,$g$-fusion frame for \,$H$.
\item[$(II)$]If for some \,$U \,\in\, G\,\mathcal{B}(\,H\,)$, \,$\left(\,\Lambda,\,  U\,\Lambda\,\right)$\, is a bi-$g$-fusion frame for \,$H$, then \,$\Lambda$\, is a \,$U$-controlled \,$g$-fusion frame for \,$H$. 
\item[$(III)$]If for some \,$T,\, U \,\in\, G\,\mathcal{B}(\,H\,)$, \,$\left(\,T\,\Lambda,\,  U\,\Lambda\,\right)$\, is a bi-$g$-fusion frame for \,$H$, then \,$\Lambda$\, is a \,$(\,T,\,U\,)$-controlled \,$g$-fusion frame for \,$H$.
\end{itemize}
\end{remark}

\begin{example}
Let \,$H$\, be a separable Hilbert space and \,$\left(\,\left\{\,f_{i}\,\right\}_{i \,\in\, I},\, \left\{\,g_{i}\,\right\}_{i \,\in\, I}\,\right)$\, be a bi-frame for\,$H$\, with bounds \,$A$\,and \,$B$.\,Suppose that \,$V_{i} \,=\, \overline{Span\,f_{i}}$\, and \,$W_{i} \,=\, \overline{Span\,g_{i}}$.\,Define \,$\Lambda_{i}\,f \,=\, \left<\,f,\, f_{i}\,\right>$\, and \,$\Gamma_{i}\,f \,=\, \left<\,g_{i},\, f\,\right>$\, for all \,$f \,\in\, H$\, and \,$i \,\in\,I$.\,Then for each \,$f \,\in\, H$, we get 
\begin{align*}
\sum\limits_{\,i \,\in\, I}\, v^{\,2}_{i}\,\left<\,\Lambda_{i}\,P_{\,V_{i}}\,f,\, \Gamma_{i}\,P_{\,W_{i}}\,f\,\right> \,=\, \sum\limits_{\,i \,\in\, I}\, v^{\,2}_{i}\,\left<\,f,\, f_{i}\,\right>\,\left<\,g_{i}\,,\, f\,\right>. 
\end{align*} 
This shows that  \,$\left(\,\Lambda,\, \Gamma\,\right) \,=\, \left(\,\left\{\,\left(\,V_{i},\, \Lambda_{i},\, v_{i}\,\right)\,\right\}_{i \,\in\, I}\,,\,\,\left\{\,\left(\,W_{i},\, \Gamma_{i},\, v_{i}\,\right)\,\right\}_{i \,\in\, I}\,\right)$\, is a bi-$g$-fusion frame for \,$H$\, with bounds \,$A$\,and \,$B$.
\end{example}

\begin{definition}
The bi-$g$-fusion frame operator \,$S_{\Lambda,\, \Gamma} \,:\, H \,\to\, H$\; is defined as follows:
\[S_{\Lambda,\, \Gamma}\,(\,f\,) \,=\, \sum\limits_{\,i \,\in\, I}\, v_{i}^{\,2}\; P_{\,W_{i}}\, \Gamma_{i}^{\,\ast}\; \Lambda_{i}\, P_{\,V_{i}}\,(\,f\,)\; \;\forall\; f \,\in\, H.\]
\end{definition}

\begin{theorem}
Let \,$\left(\,\Lambda,\, \Gamma\,\right)$\, be a bi-$g$-fusion frame for \,$H$\, with bounds \,$A$\, and \,$B$.\,Then the bi-$g$-fusion frame operator \,$S_{\Lambda,\, \Gamma}$\, is well-defined, bounded, positive, invertible and moreover \,$S^{\,\ast}_{\Lambda,\, \Gamma} \,=\, S_{\Gamma,\,\Lambda}$. 
\end{theorem}

\begin{proof}
For each \,$f,\; g \,\in\, H$, we have
\begin{align}
\left<\,S_{\Lambda,\,\Gamma\,}\,(\,f\,) \,,\, g\,\right> &\,=\, \left<\,\sum\limits_{\,i \,\in\, I}\, v^{2}_{i}\, P_{\,W_{i}}\, \Gamma_{i}^{\,\ast}\; \Lambda_{i}\, P_{\,V_{i}}\,(\,f\,) \;,\; g\,\right>\nonumber\\
& =\; \sum\limits_{\,i \,\in\, I}\, v^{2}_{i}\,\left<\,\Lambda_{i}\, P_{\,V_{i}}\,(\,f\,) \;,\; \Gamma_{i}\,P_{\,W_{i}}\,(\,g\,)\,\right>.\label{eq3}
\end{align} 
By the Cauchy-Schwarz inequality, we obtain 
\begin{align}
&\left|\,\left<\,S_{\Lambda,\,\Gamma}\,(\,f\,) \,,\, g\,\right>\,\right|\nonumber\\
& \,\leq\, \left(\,\sum\limits_{\,i \,\in\, I}\, v_{\,i}^{\,2}\; \left\|\,\Gamma_{i}\, P_{\,W_{i}}\,(\,g\,) \,\right\|^{\,2}\,\right)^{\,1\,/\, 2}\; \left(\,\sum\limits_{\,i \,\in\, I}\, v_{\,i}^{\,2}\;  \left\|\,\Lambda_{i}\, P_{\,V_{i}}\,(\,f\,) \,\right\|^{\,2}\,\right)^{\,1\,/\, 2}\label{eq4}
\end{align}
Now, for each \,$f \,\in\, H$, we have
\begin{align*}
&\left\|\,S_{\Lambda,\,\Gamma}\,f\,\right\| \;=\; \sup\left\{\,\left|\,\left<\,S_{\Lambda,\,\Gamma\,}(\,f\,) \,,\, f\,\right>\,\right| \;:\; \|\,f\,\| \;=\; 1\,\right\}\\
&=\,\sup\left\{\,\left|\,\sum\limits_{\,i \,\in\, I}\, v^{2}_{i}\,\left<\,\Lambda_{i}\, P_{\,V_{i}}\,(\,f\,) \;,\; \Gamma_{i}\,P_{\,W_{i}}\,(\,f\,)\,\right>\,\right| \;:\; \|\,f\,\| \;=\; 1\,\right\} \,\leq\, B.
\end{align*} 
Also, for each \,$f,\, g \,\in\, H$, we have
\begin{align*}
\left<\,S_{\Gamma,\,\Lambda}\,(\,f\,) \;,\; g\,\right>& \,=\, \left<\,\sum\limits_{\,i \,\in\, I}\, v^{2}_{i}\, P_{\,V_{i}}\, \Lambda_{i}^{\,\ast}\; \Gamma_{i}\, P_{\,W_{i}}\,(\,f\,) \;,\; g\,\right>\\
&=\,\sum\limits_{\,i \,\in\, I}\, v^{2}_{\,i}\, \left<\,f \;,\;  P_{\,W_{i}}\,\Gamma_{i}^{\,\ast}\; \Lambda_{i}\, P_{\,V_{i}}\,(\,g\,) \,\right>\\
& \,=\, \left<\,f \;,\;  \sum\limits_{\,i \,\in\, I}\, v^{2}_{\,i}\,P_{\,W_{i}}\,\Gamma_{i}^{\,\ast}\; \Lambda_{i}\, P_{\,V_{i}}\,(\,g\,) \,\right>\\
&\,=\, \left<\,f \;,\; S_{\Lambda,\,\Gamma}\,(\,g\,)\,\right>.
\end{align*}
Thus, \,$S_{\Gamma,\,\Lambda}^{\,\ast} \,=\, S_{\Lambda,\,\Gamma}$.\,Since \,$\left(\,\Lambda,\, \Gamma\,\right)$\, is a bi-$g$-fusion frame for \,$H$\, with bounds \,$A$\, and \,$B$, for each \,$f,\, g \,\in\, H$, we get \,$A\,\left<\,f,\, f\,\right>\,\leq\, \left<\,S_{\Lambda,\, \Gamma}\,f,\; f\,\right> \,\leq\,B\,\left<\,f,\, f\,\right>$.\,Therefore, \,$A\,I_{H}\,\leq\, S_{\Lambda,\, \Gamma}\,\leq\,B\,I_{H}$.\,This shows that \,$S_{\Lambda,\, \Gamma}$\, is a positive operator.\,Furthermore, it is easy to verify that \,$S_{\Lambda,\, \Gamma}$\, is invertible operator.\,This completes the proof.  
\end{proof}

\begin{theorem}
The pair \,$\left(\,\Lambda,\, \Gamma\,\right)$\, is a bi-$g$-fusion frame for \,$H$\, if and only if \,$\left(\,\Gamma,\, \Lambda\,\right)$\, is a bi-$g$-fusion frame for \,$H$.
\end{theorem}

\begin{proof}
Let \,$\left(\,\Lambda,\, \Gamma\,\right)$\, is a bi-$g$-fusion frame for \,$H$\,with bounds \,$A$\, and \,$B$.\,Then for each \,$f \,\in\, H$, we have   
\[A\,\|\,f\,\|^{\,2} \,\leq\, \sum\limits_{\,i \,\in\, I}\, v^{\,2}_{i}\,\left<\,\Lambda_{i}\,P_{\,V_{i}}\,f,\, \Gamma_{i}\,P_{\,W_{i}}\,f\,\right> \,\leq\, \,B\,\|\,f \,\|^{\,2}\]
Now, for each \,$f \,\in\, H$, we can write
\begin{align*}
&\sum\limits_{\,i \,\in\, I}\, v^{\,2}_{i}\,\left<\,\Lambda_{i}\,P_{\,V_{i}}\,f,\, \Gamma_{i}\,P_{\,W_{i}}\,f\,\right> \,=\, \overline{\,\sum\limits_{\,i \,\in\, I}\, v^{\,2}_{i}\,\left<\,\Lambda_{i}\,P_{\,V_{i}}\,f,\,  \Gamma_{i}\,P_{\,W_{i}}\,f\,\right>}\\
&=\,\sum\limits_{\,i \,\in\, I}\, v^{\,2}_{i}\,\overline{\,\left<\,\Lambda_{i}\,P_{\,V_{i}}\,f,\, \Gamma_{i}\,P_{\,W_{i}}\,f\,\right>} \,=\, \sum\limits_{\,i \,\in\, I}\, v^{\,2}_{i}\,\left<\,\Gamma_{i}\,P_{\,W_{i}}\,f,\, \Lambda_{i}\,P_{\,V_{i}}\,f\,\right>. 
\end{align*}
This shows that \,$\left(\,\Gamma,\, \Lambda\,\right)$\, is a bi-$g$-fusion frame for \,$H$\, with bounds \,$A$\, and \,$B$.

Similarly, we can prove the converse part of this Theorem.
\end{proof}

\begin{remark}
The pair \,$\left(\,\Lambda,\, \Gamma\,\right)$\, is a bi-$g$-fusion frame for \,$H$\, if and only if the bi-$g$-fusion frame operator \,$S_{\Lambda,\, \Gamma}$\, is well-defined, bounded and surjective. 
\end{remark}

\begin{remark}
Let \,$\left(\,\Lambda,\, \Gamma\,\right)$\, be a bi-$g$-fusion frame for \,$H$\, with bi-$g$-fusion frame operator \,$S_{\Lambda,\, \Gamma}$.\,Then for any \,$f \,\in\, H$, we get
\begin{align*}
 f &\,=\, S_{\Lambda,\, \Gamma}\,S^{\,-\, 1}_{\Lambda,\, \Gamma}\,f \,=\, \sum\limits_{\,i \,\in\, I}\, v_{i}^{\,2}\; P_{\,W_{i}}\, \Gamma_{i}^{\,\ast}\; \Lambda_{i}\, P_{\,V_{i}}\,S^{\,-\, 1}_{\Lambda,\, \Gamma}\,(\,f\,)\\
 &=\,S^{\,-\, 1}_{\Lambda,\, \Gamma}\,S_{\Lambda,\, \Gamma}\,f \,=\,  \,=\, \sum\limits_{\,i \,\in\, I}\,v^{\,2}_{\,i}\,S^{\,-\, 1}_{\Lambda,\, \Gamma}\,P_{\,W_{i}}\,\Gamma^{\,\ast}_{i}\,\Lambda_{i}\,P_{\,V_{i}}\,(\,f\,).
\end{align*}  
This is called reconstructions formula for bi-$g$-fusion frame.
\end{remark}

\begin{theorem}
Let \,$\left(\,\Lambda,\, \Gamma\,\right)$\, be a bi-$g$-fusion Bessel sequence in \,$H$\, with bi-$g$-fusion frame operator \,$S_{\Lambda,\, \Gamma}$.\,Then \,$\left(\,\Lambda,\, \Gamma\,\right)$\, is a bi-$g$-fusion frame for \,$H$\, if and only if there exists \,$\alpha \,>\, 0$\,such that \,$S_{\Lambda,\, \Gamma}\,\geq\, \alpha\,I_{H}$.
\end{theorem}

\begin{proof}
Let \,$\left(\,\Lambda,\, \Gamma\,\right)$\, is a bi-$g$-fusion frame for \,$H$\,with bounds \,$A$\, and \,$B$.\,Then for each \,$f \,\in\, H$, we have 
\begin{align*}
&A\,\|\,f\,\|^{\,2} \,\leq\, \sum\limits_{\,i \,\in\, I}\, v^{\,2}_{i}\,\left<\,\Lambda_{i}\,P_{\,V_{i}}\,f,\, \Gamma_{i}\,P_{\,W_{i}}\,f\,\right> \,\leq\, \,B\,\|\,f \,\|^{\,2}\\
&\Rightarrow\, A\,\|\,f\,\|^{\,2} \,\leq\, \left<\,\sum\limits_{\,i \,\in\, I}\, v^{2}_{i}\, P_{\,W_{i}}\, \Gamma_{i}^{\,\ast}\; \Lambda_{i}\, P_{\,V_{i}}\,(\,f\,) \;,\; f\,\right>\, \leq\,B\,\|\,f\,\|^{\,2} \\
&\Rightarrow A\,\left<\,f,\, f\,\right>\,\leq\, \left<\,S_{\Lambda,\, \Gamma}\,f,\; f\,\right> \,\leq\,B\,\left<\,f,\, f\,\right>.
\end{align*}  
This shows that \,$\alpha\,\,I_{H}\,\leq\, S_{\Lambda,\, \Gamma} \,$,  where \,$\alpha\,=\, A$.

Converse part is obvious.   
\end{proof}

\begin{theorem}
Let \,$\left(\,\Lambda,\,  \Gamma\,\right)$\, be a bi-$g$-fusion Bessel sequence in \,$H$\,with bound \,$B$\, and \,$S_{\Lambda,\, \Gamma}$\; be the bi-$g$-fusion frame operator.\,Then the following statements are equivalent:
\begin{description}
\item[$(i)$]\; $S_{\Lambda,\, \Gamma}$\; is bounded below.
\item[$(ii)$] There exists \,$K \,\in\, \mathcal{B}\,(\,H\,)$\; such that \;$\left\{\,T_{i}\,\right\}_{i \,\in\, I}$\; is a resolution of the identity operator on \,$H$, where \;$T_{i} \,=\, v^{\,2}_{i}\, K\, P_{\,W_{i}}\, \Gamma_{i}^{\,\ast}\; \Lambda_{i}\, P_{\,V_{i}} \;,\; i \,\in\, I$.  
\end{description} 
\end{theorem}

\begin{proof}
$(\,i\,) \,\Rightarrow\, (\,ii\,)$\; Suppose that \,$S_{\Lambda,\, \Gamma}$\, is bounded below.\;Then for each \,$f \,\in\, H$, there exists \,$A \,>\, 0$\; such that 
\begin{align*}
\|\,f\,\|^{\,2} \,\leq\, A\;\left\|\,S_{\Lambda,\, \Gamma}\,f\,\right\|^{\,2} &\,\Rightarrow\, \left<\,I_{H}\,f \,,\, f\,\right> \,\leq\, A\, \left<\,S_{\Lambda,\, \Gamma}^{\,\ast}\, S_{\Lambda,\, \Gamma}\,f \,,\, f\,\right>\\
& \,\Rightarrow\,I^{\,\ast}_{H}\,I_{\,H} \,\leq\, A\, S_{\Lambda,\, \Gamma}^{\,\ast}\,S_{\Lambda,\, \Gamma}.
\end{align*} 
So, by Theorem \ref{th1}, there exists \,$K \,\in\, \mathcal{B}\,(\,H\,)$\, such that \,$K\, S_{\Lambda,\, \Gamma} \,=\, I_{\,H}$.\\Therefore, for each \,$f \,\in\, H$, we have
\begin{align*}
 f &\,=\, K\, S_{\Lambda,\, \Gamma}\,(\,f\,) \,=\, K\, \left(\,\sum\limits_{\,i \,\in\, I}\, v^{\,2}_{i}\, P_{\,W_{i}}\, \Gamma_{i}^{\,\ast}\; \Lambda_{i}\, P_{\,V_{i}}\,(\,f\,)\,\right)\\
&\,=\, \sum\limits_{\,i \,\in\, I}\, v^{\,2}_{i}\, K\, P_{\,W_{i}}\, \Gamma_{i}^{\,\ast}\; \Lambda_{i}\, P_{\,V_{i}} \,=\, \sum\limits_{\,i \,\in\, I}\, T_{i}\,(\,f\,)
\end{align*}
 
 and hence \,$\left\{\,T_{\,j}\,\right\}_{i \,\in\, I}$\; is a resolution of the identity operator on \,$H$, where \,$T_{i} \,=\, v^{\,2}_{i}\, K\, P_{\,W_{i}}\, \Gamma_{i}^{\,\ast}\; \Lambda_{i}\, P_{\,V_{i}}$.\\\\
$(\,ii\,) \,\Rightarrow\, (\,i\,)$\; Since \;$\left\{\,T_{i}\,\right\}_{i \,\in\, I}$\; is a resolution of the identity operator on \,$H$, for any \,$f \,\in\, H$, we have
\begin{align*}
f &\,=\, \sum\limits_{\,i \,\in\, I}\, T_{\,i}\,(\,f\,) \,=\, \sum\limits_{\,i \,\in\, I}\, v^{\,2}_{i}\, K\, P_{\,W_{i}}\, \Gamma_{i}^{\,\ast}\; \Lambda_{i}\, P_{\,V_{i}}\,(\,f\,)\\
& \,=\, K\,\left(\,\sum\limits_{\,i \,\in\, I}\, v^{\,2}_{i}\, P_{\,W_{i}}\, \Gamma_{i}^{\,\ast}\; \Lambda_{i}\, P_{\,V_{i}}\,(\,f\,)\,\right) \,=\, K\, S_{\Lambda,\, \Gamma}\,(\,f\,).
\end{align*}
Thus, \,$I_{H} \,=\, K\,S_{\Lambda,\, \Gamma}$.\;So, by Theorem \ref{th1}, there exists some \,$\alpha \,>\, 0$\; such that \,$I_{H}\, I_{H}^{\,\ast} \,\leq\; \alpha\; S_{\Lambda,\, \Gamma}\, S_{\Lambda,\, \Gamma}^{\,\ast}$\; and hence \;$S_{\Lambda,\, \Gamma}$\; is bounded below.\,This completes the proof.
\end{proof}

\begin{theorem}\label{thm3.13}
Let \,$U \,\in\, \mathcal{B}(\,H\,)$\, be an invertible operator and \,$\left(\,\Lambda,\, \Gamma\,\right)$\, be a bi-$g$-fusion frame for \,$H$.\,Then 
\[\left(\,\left\{\,\left(\,U\,V_{i},\, \Lambda_{i}\,P_{\,V_{i}}\,U^{\,\ast},\, v_{i}\,\right)\,\right\}_{i \,\in\, I}\,,\, \,\left\{\,\left(\,U\,W_{i},\, \Gamma_{i}\,P_{\,W_{i}}\,U^{\,\ast},\, v_{i}\,\right)\,\right\}_{i \,\in\, I}\,\right)\] is a bi-$g$-fusion frame for \,$H$.  
\end{theorem}

\begin{proof}
Let \,$\left(\,\Lambda,\, \Gamma\,\right)$\, be a bi-$g$-fusion frame for \,$H$\, with bounds \,$A$\, and \,$B$.\,Since \,$U$\, is invertible, for each \,$f \,\in\, H$, we get
\[\|\,f \,\|^{\,2}\,=\, \left\|\,\left(\,U^{\,-\,1}\,\right)^{\,\ast}\,U^{\,\ast}\,f \,\right\|^{\,2} \,\leq\, \left\|\,U^{\,-\,1} \,\right\|^{\,2}\,\|\,f \,\|^{\,2}.\]
Now, for each \,$f \,\in\, H$, we have 
\begin{align}
&\sum\limits_{\,i \,\in\, I}\, v^{\,2}_{i}\,\left<\,\Lambda_{i}\,P_{\,V_{i}}\,U^{\,\ast}\,P_{\,U\,V_{i}}\,f\,,\, \Gamma_{i}\,P_{\,W_{i}}\,U^{\,\ast}\,P_{\,U\,W_{i}}\,f\,\right>\nonumber\\
&=\,\sum\limits_{\,i \,\in\, I}\, v^{\,2}_{i}\,\left<\,\Lambda_{i}\,P_{\,V_{i}}\,U^{\,\ast}\,f\,,\, \Gamma_{i}\,P_{\,W_{i}}\,U^{\,\ast}\,f\,\right>\label{eqn3.13} \\
&\leq\,B\,\left\|\,U^{\,\ast}\,f \,\right\|^{\,2} \,\leq\, B\,\left\|\,U \,\right\|^{\,2}\,\|\,f \,\|^{\,2}.\nonumber
\end{align}
On the other hand, from (\ref{eqn3.13}), for each \,$f \,\in\, H$, we get
\begin{align*}
&\sum\limits_{\,i \,\in\, I}\, v^{\,2}_{i}\,\left<\,\Lambda_{i}\,P_{\,V_{i}}\,U^{\,\ast}\,P_{\,U\,V_{i}}\,f\,,\, \Gamma_{i}\,P_{\,W_{i}}\,U^{\,\ast}\,P_{\,U\,W_{i}}\,f\,\right>\\
&=\,\sum\limits_{\,i \,\in\, I}\, v^{\,2}_{i}\,\left<\,\Lambda_{i}\,P_{\,V_{i}}\,U^{\,\ast}\,f\,,\, \Gamma_{i}\,P_{\,W_{i}}\,U^{\,\ast}\,f\,\right>\\
&\geq\, A\,\left\|\,U^{\,\ast}\,f \,\right\|^{\,2} \,\geq\, A\,\left\|\,U^{\,-\,1} \,\right\|^{\,-\,2}\,\|\,f \,\|^{\,2}.
\end{align*} 
This completes the proof.
\end{proof}

\begin{corollary}
Let \,$\left(\,\Lambda,\, \Gamma\,\right)$\, be a bi-$g$-fusion frame for \,$H$\, with the corresponding bi-$g$-fusion frame operator \,$S_{\Lambda,\, \Gamma}$.\,Then the family \,$\left(\,\overline{\,\Lambda},\, \overline{\,\Gamma}\,\right)$
\[=\, \left(\,\left\{\,\left(\,S^{\,-\, 1}_{\Lambda,\, \Gamma}\,V_{i},\, \Lambda_{i}\,P_{\,V_{i}}\,S^{\,-\, 1}_{\Lambda,\, \Gamma},\, v_{i}\,\right)\,\right\}_{i \,\in\, I}\,,\, \,\left\{\,\left(\,S^{\,-\, 1}_{\Lambda,\, \Gamma}\,W_{i},\, \Gamma_{i}\,P_{\,W_{i}}\,S^{\,-\, 1}_{\Lambda,\, \Gamma},\, v_{i}\,\right)\,\right\}_{i \,\in\, I}\,\right)\] is a bi-$g$-fusion frame for \,$H$.  
\end{corollary}

The bi-$g$-fusion frame \,$\left(\,\overline{\,\Lambda},\, \overline{\,\Gamma}\,\right)$\,is called the canonical dual bi \,$g$-fusion frame of \,$\left(\,\Lambda,\, \Gamma\,\right)$.

%=====================================
\section{Bi-$g$-fusion frames in \,$H \,\otimes\, K$\,}
%=====================================

In this section, \;$H \;\text{and}\; K$\, are considered to be separable Hilbert spaces with associated inner products \,$\left <\,\cdot \,,\, \cdot\,\right>_{1} \;\text{and}\; \left <\,\cdot \,,\, \cdot\,\right>_{2}$.\,$\left\{\,W_{j}\,\right\}_{ j \,\in\, J}$\, is the collections of closed subspaces of \,$K$, where \,$J$\, is subset of  integers \,$\mathbb{Z}$.\,$\left\{\,K_{j}\,\right\}_{ j \,\in\, J}$\, is the collections of Hilbert spaces.\,$\left\{\,\Gamma_{j} \,\in\, \mathcal{B}\,(\,K \,,\, K_{j}\,)\,\right\}_{j \,\in\, J}$\, denotes the sequence of operators.\,As an above, we can define the space \,$l^{\,2}\left(\,\left\{\,K_{j}\,\right\}_{ j \,\in\, J}\,\right)$.\,Now, we present bi-$g$-fusion frame in tensor product of Hilbert spaces \,$H \,\otimes\, K$. 

\begin{definition}
Let \,$\left\{\,v_{\,i}\,\right\}_{\, i \,\in\, I},\;\left\{\,w_{\,j}\,\right\}_{ j \,\in\, J}$\, be two families of positive weights i\,.\,e., \,$v_{\,i} \,>\, 0\, \;\forall\; i \,\in\, I,\; \,w_{\,j} \,>\, 0\, \;\forall\; j \,\in\, J$.\,Suppose \,$\Lambda_{i},\, \Lambda_{i}^{\,\prime} \,\in\, \mathcal{B}\,(\,H,\, H_{i}\,)\,$\, and \,$\Gamma_{j},\, \Gamma_{j}^{\,\prime} \,\in\, \mathcal{B}\,(\,K,\,K_{j}\,)$\, for each \,$i \,\in\, I$, \,$j \,\in\, J$.\,Then the pairs of families \,$\left\{\,\left(\,V_{i} \,\otimes\, W_{j},\, \Lambda_{i} \,\otimes\, \Gamma_{j},\, v_{i}\,w_{j}\,\right)\,\right\}_{\,i,\,j}$\, and \,$\left\{\,\left(\,V_{i}^{\,\prime} \,\otimes\, W_{j}^{\,\prime},\, \Lambda_{i}^{\,\prime} \,\otimes\, \Gamma_{j}^{\,\prime},\, v_{i}\,w_{j}\,\right)\,\right\}_{\,i,\,j}$\, is said to be a bi-$g$-fusion frame for \,$H \,\otimes\, K$\, with respect to \,$\left\{\,H_{i} \,\otimes\, K_{j}\,\right\}_{\,i,\,j}$\, if there exist constants \,$0 \,<\, A \,\leq\, B \,<\, \infty$\, such that
\begin{align}
&A\, \left\|\,f \,\otimes\, g\,\right\|^{\,2}\nonumber\\
&\leq\,\sum\limits_{i,\, j}\,v^{\,2}_{\,i}\,w^{\,2}_{\,j}\,\left<\left(\,\Lambda_{i} \,\otimes\, \Gamma_{j}\,\right)\,P_{\,V_{\,i} \,\otimes\, W_{\,j}}\,(\,f \,\otimes\, g\,),\, \left(\,\Lambda_{i}^{\,\prime} \,\otimes\, \Gamma_{j}^{\,\prime}\,\right)\,P_{\,V_{\,i}^{\,\prime} \,\otimes\, W_{\,j}^{\,\prime}}\,(\,f \,\otimes\, g\,)\,\right>\nonumber\\
&\leq\,B\, \left\|\,f \,\otimes\, g\,\right\|^{\,2}\; \;\forall\; f \,\otimes\, g \,\in\, H \,\otimes\, K \label{eqn2.1}
\end{align}
where \,$P_{\,V_{i} \,\otimes\, W_{j}}$\, and \,$P_{\,V_{\,i}^{\,\prime} \,\otimes\, W_{\,j}^{\,\prime}}$\, are  the orthogonal projections of \,$H \,\otimes\, K$\, onto \,$V_{\,i} \,\otimes\, W_{\,j}$\, and \,$V_{\,i}^{\,\prime} \,\otimes\, W_{\,j}^{\,\prime}$, respectively.\,The constants \,$A$\, and \,$B$\, are called the frame bounds of \,$\Lambda \,\otimes\, \Gamma$.\;If \,$A \,=\, B$\, then it is called a tight bi-$g$-fusion frame.\,If the family satisfies the right inequality of (\ref{eqn2.1}) then it is called a bi-$g$-fusion Bessel sequence in \,$H \,\otimes\, K$\, with bound \,$B$.   
\end{definition}

\begin{remark}
For \,$i \,\in\, I$\, and \,$j \,\in\, J$, define the space \,$l^{\,2}\,\left(\,\left\{\,H_{i} \,\otimes\, K_{j}\,\right\}\,\right)$
\[\,=\, \left\{\,\left\{\,f_{\,i} \,\otimes\, g_{\,j}\,\right\} \,:\, f_{\,i} \,\otimes\, g_{\,j} \,\in\, H_{i} \,\otimes\, K_{j}, \;\&\; \;\sum\limits_{i,\, j}\,\left\|\,f_{\,i} \,\otimes\, g_{\,j}\,\right\|^{\,2} \,<\, \infty\,\right\}\]
with the inner product given by 
\begin{align*}
&\left<\,\left\{\,f_{\,i} \otimes g_{\,j}\,\right\} \,,\, \left\{\,f^{\,\prime}_{\,i} \otimes g^{\,\prime}_{\,j}\,\right\}\,\right>_{l^{\,2}} \,=\, \sum\limits_{i,\, j}\,\left<\,f_{\,i} \otimes g_{\,j}\, \,,\, f^{\,\prime}_{\,i} \otimes g^{\,\prime}_{\,j}\,\right>\\
&=\,\sum\limits_{i,\, j}\,\left<\,f_{\,i} \,,\, f^{\,\prime}_{\,i}\,\right>_{H_{i}}\,\left<\,g_{\,j} \,,\, g^{\,\prime}_{\,j}\,\right>_{K_{j}} \,=\, \left(\,\sum\limits_{\,i \,\in\, I}\,\left<\,f_{\,i} \,,\, f^{\,\prime}_{\,i}\,\right>_{H_{i}}\,\right)\,\left(\,\sum\limits_{\,j \,\in\, J}\,\left<\,g_{\,j} \,,\, g^{\,\prime}_{\,j}\,\right>_{K_{j}}\,\right)\\
&=\, \left<\,\left\{\,f_{\,i}\,\right\}_{ i \,\in\, I} \,,\, \left\{\,f^{\,\prime}_{\,i}\,\right\}_{ i \,\in\, I}\,\right>_{l^{\,2}\left(\,\left\{\,H_{i}\,\right\}_{i \,\in\, I}\,\right)}\,\left<\,\left\{\,g_{\,j}\,\right\}_{ j \,\in\, J} \,,\, \left\{\,g^{\,\prime}_{\,j}\,\right\}_{ j \,\in\, J}\,\right>_{l^{\,2}\left(\,\left\{\,K_{j}\,\right\}_{ j \,\in\, J}\,\right)}.
\end{align*}
The space \,$l^{\,2}\,\left(\,\left\{\,H_{i} \,\otimes\, K_{j}\,\right\}\,\right)$\, is complete with the above inner product.\;Then it becomes a Hilbert space with respect to the above inner product.  
\end{remark}

\begin{remark}
Since \,$\left\{\,V_{i}\,\right\}_{ i \,\in\, I},\; \left\{\,W_{j}\,\right\}_{ j \,\in\, J}$\, and \,$\left\{\,V_{i} \,\otimes\, W_{j} \right\}_{i,\,j}$\, are the families of closed subspaces of the Hilbert spaces \,$H,\, K$\, and \,$H \,\otimes\, K$, respectively, it is easy to verify that \,$P_{\,V_{i} \,\otimes\, W_{j}} \,=\, P_{\,V_{i}} \,\otimes\, P_{\,W_{j}}$.
\end{remark}

\begin{theorem}\label{th1.2}
The pair of families \,$\left(\,\Lambda,\, \Lambda^{\,\prime}\right)$\, and \,$\left(\,\Gamma,\,\Gamma^{\,\prime}\,\right)$\, are bi-$g$-fusion frames for \,$H$\, and \,$K$\, with respect to \,$\left\{\,H_{i}\,\right\}_{\, i \,\in\, I}$\, and \,$\left\{\,K_{j}\,\right\}_{\, j \,\in\, J}$, respectively if and only if the family \,$\left(\,\Lambda \,\otimes\, \Gamma,\, \Lambda^{\,\prime} \,\otimes\, \Gamma^{\,\prime}\,\right)$\, is a bi $g$-fusion frame for \,$H \,\otimes\, K$\, with respect to \,$\left\{\,H_{i} \,\otimes\, K_{j}\,\right\}_{\,i,\,j}$. 
\end{theorem}

\begin{proof}
First we suppose that \,$\left(\,\Lambda \,\otimes\, \Gamma,\, \Lambda^{\,\prime} \,\otimes\, \Gamma^{\,\prime}\,\right)$\, is a bi $g$-fusion frame for \,$H \,\otimes\, K$\, with respect to \,$\left\{\,H_{i} \,\otimes\, K_{j}\,\right\}_{\,i,\,j}$\, having bounds \,$A$\, and \,$B$.\,Then for each \,$f \,\otimes\, g \,\in\, H \,\otimes\, K \,-\, \{\,\theta \,\otimes\, \theta\,\}$, using (\ref{eq1.001}) and (\ref{eq1.0001}), we have
\begin{align*}
&A\, \left\|\,f \,\otimes\, g\,\right\|^{\,2}\\
&\leq\,\sum\limits_{i,\, j}\,v^{\,2}_{\,i}\,w^{\,2}_{\,j}\,\left<\left(\,\Lambda_{i} \,\otimes\, \Gamma_{j}\,\right)\,P_{\,V_{\,i} \,\otimes\, W_{\,j}}\,(\,f \,\otimes\, g\,),\, \left(\,\Lambda_{i}^{\,\prime} \,\otimes\, \Gamma_{j}^{\,\prime}\,\right)\,P_{\,V_{\,i}^{\,\prime} \,\otimes\, W_{\,j}^{\,\prime}}\,(\,f \,\otimes\, g\,)\,\right>\\
&\hspace{2cm}\leq\,B\, \left\|\,f \,\otimes\, g\,\right\|^{\,2}.\\
&\Rightarrow\, A\, \left\|\,f \,\otimes\, g\,\right\|^{\,2}\leq\,\sum\limits_{i,\, j}\,v^{\,2}_{\,i}\,w^{\,2}_{\,j}\,\left<\,\Lambda_{i}\, P_{\,V_{\,i}}\,f \,\otimes\, \Gamma_{j}\,P_{\,W_{\,j}}\,g\,,\, \Lambda_{i}^{\,\prime}\,P_{\,V_{\,i}^{\,\prime}}\,f \,\otimes\, \Gamma_{j}^{\,\prime}\,P_{\,W_{\,j}^{\,\prime}}\,g\,\right>\\
&\hspace{2cm}\leq\,B\, \left\|\,f \,\otimes\, g\,\right\|^{\,2}.\\
&\Rightarrow\, A\,\|\,f\,\|_{1}^{\,2}\,\|\,g\,\|_{2}^{\,2}\leq\,\sum\limits_{i,\, j}\,v^{\,2}_{\,i}\,w^{\,2}_{\,j}\,\left<\,\Lambda_{i}\, P_{\,V_{\,i}}\,f,\, \Lambda_{i}^{\,\prime}\,P_{\,V_{\,i}^{\,\prime}}\,f\,\right>_{1}\,\left<\,\Gamma_{j}\,P_{\,W_{\,j}}\,g,\, \Gamma_{j}^{\,\prime}\,P_{\,W_{\,j}^{\,\prime}}\,g\,\right>_{2}\\
&\hspace{2cm}\leq\,B\,\|\,f\,\|_{1}^{\,2}\,\|\,g\,\|_{2}^{\,2}.\\
&\Rightarrow\, A\,\|\,f\,\|_{1}^{\,2}\,\|\,g\,\|_{2}^{\,2}\leq\,\sum\limits_{i \,\in\, I}\,v^{\,2}_{\,i}\,\left<\,\Lambda_{i}\, P_{\,V_{\,i}}\,f,\, \Lambda_{i}^{\,\prime}\,P_{\,V_{\,i}^{\,\prime}}\,f\,\right>_{1}\,\sum\limits_{j \,\in\, J}\,w^{\,2}_{\,j}\,\left<\,\Gamma_{j}\,P_{\,W_{\,j}}\,g,\, \Gamma_{j}^{\,\prime}\,P_{\,W_{\,j}^{\,\prime}}\,g\,\right>_{2}\\
&\hspace{2cm}\leq\,B\,\|\,f\,\|_{1}^{\,2}\,\|\,g\,\|_{2}^{\,2}.\\
\end{align*}
Since \,$f \,\otimes\, g$\, is non-zero vector, \,$f$\, and \,$g$\, are also non-zero vectors and therefore 
\[\sum\limits_{i \,\in\, I}\,v^{\,2}_{\,i}\,\left<\,\Lambda_{i}\, P_{\,V_{\,i}}\,f,\, \Lambda_{i}^{\,\prime}\,P_{\,V_{\,i}^{\,\prime}}\,f\,\right>_{1},\, \,\sum\limits_{j \,\in\, J}\,w^{\,2}_{\,j}\,\left<\,\Gamma_{j}\,P_{\,W_{\,j}}\,g,\, \Gamma_{j}^{\,\prime}\,P_{\,W_{\,j}^{\,\prime}}\,g\,\right>_{2}\] are non-zero.\,Therefore from the above inequality we get
\begin{align*}
&\dfrac{A\,\left\|\,g \,\right\|_{\,2}^{\,2}\,\left\|\,f \,\right\|_{\,1}^{\,2}}{\sum\limits_{j \,\in\, J}\,w^{\,2}_{\,j}\,\left<\,\Gamma_{j}\,P_{\,W_{\,j}}\,g,\, \Gamma_{j}^{\,\prime}\,P_{\,W_{\,j}^{\,\prime}}\,g\,\right>_{2}} \,\leq\, \sum\limits_{i \,\in\, I}\,v^{\,2}_{\,i}\,\left<\,\Lambda_{i}\, P_{\,V_{\,i}}\,f,\, \Lambda_{i}^{\,\prime}\,P_{\,V_{\,i}^{\,\prime}}\,f\,\right>_{1}\\
&\hspace{4cm} \,\leq\, \dfrac{B\,\left\|\,g \,\right\|_{\,2}^{\,2}\,\left\|\,f \,\right\|_{\,1}^{\,2}}{\sum\limits_{j \,\in\, J}\,w^{\,2}_{\,j}\,\left<\,\Gamma_{j}\,P_{\,W_{\,j}}\,g,\, \Gamma_{j}^{\,\prime}\,P_{\,W_{\,j}^{\,\prime}}\,g\,\right>_{2}}\\
&\Rightarrow\, A_{1} \,\left\|\,f \,\right\|_{\,1}^{\,2} \,\leq\, \sum\limits_{i \,\in\, I}\,v^{\,2}_{\,i}\,\left<\,\Lambda_{i}\, P_{\,V_{\,i}}\,f,\, \Lambda_{i}^{\,\prime}\,P_{\,V_{\,i}^{\,\prime}}\,f\,\right>_{1} \,\leq\, B_{1}\, \left\|\, f \, \right\|_{\,1}^{\,2}\; \;\forall\; f \,\in\, H,
\end{align*}
where 
\[A_{1} \,=\, \inf\limits_{g \,\in\, K}\,\dfrac{A\,\left\|\,g \,\right\|_{\,2}^{\,2}}{\sum\limits_{j \,\in\, J}\,w^{\,2}_{\,j}\,\left<\,\Gamma_{j}\,P_{\,W_{\,j}}\,g,\, \Gamma_{j}^{\,\prime}\,P_{\,W_{\,j}^{\,\prime}}\,g\,\right>_{2}},\,\text{and}\]
\[B_{1} \,=\, \sup\limits_{g \,\in\, K}\,\dfrac{B\,\left\|\,g \,\right\|_{\,2}^{\,2}}{\sum\limits_{j \,\in\, J}\,w^{\,2}_{\,j}\,\left<\,\Gamma_{j}\,P_{\,W_{\,j}}\,g,\, \Gamma_{j}^{\,\prime}\,P_{\,W_{\,j}^{\,\prime}}\,g\,\right>_{2}}.\]This shows that the pair \,$\left(\,\Lambda,\, \Lambda^{\,\prime}\right)$\, is a bi-$g$-fusion frame for \,$H$\, with respect to \,$\left\{\,H_{i}\,\right\}_{\, i \,\in\, I}$.\,Similarly, it can be shown that and \,$\left(\,\Gamma,\,\Gamma^{\,\prime}\,\right)$\, is bi-$g$-fusion frame for \,$K$\, with respect to \,$\left\{\,K_{j}\,\right\}_{\, j \,\in\, J}$.

Conversely, suppose that \,$\left(\,\Lambda,\, \Lambda^{\,\prime}\right)$\, and \,$\left(\,\Gamma,\,\Gamma^{\,\prime}\,\right)$\, are bi-$g$-fusion frames for \,$H$\, and \,$K$\, with bounds \,$(\,A,\, B\,)$\, and \,$(\,C,\, D\,)$, respectively.\,Then
\begin{align}
&A \,\left\|\,f \,\right\|_{\,1}^{\,2} \,\leq\, \sum\limits_{i \,\in\, I}\,v^{\,2}_{\,i}\,\left<\,\Lambda_{i}\, P_{\,V_{\,i}}\,f,\, \Lambda_{i}^{\,\prime}\,P_{\,V_{\,i}^{\,\prime}}\,f\,\right>_{1} \,\leq\, B\, \left\|\, f \, \right\|_{\,1}^{\,2}\; \;\forall\; f \,\in\, H\label{eq1}\\
&C \,\left\|\,g \,\right\|_{\,1}^{\,2} \,\leq\, \sum\limits_{j \,\in\, J}\,w^{\,2}_{\,j}\,\left<\,\Gamma_{j}\,P_{\,W_{\,j}}\,g,\, \Gamma_{j}^{\,\prime}\,P_{\,W_{\,j}^{\,\prime}}\,g\,\right>_{2} \,\leq\, D\, \left\|\, g \, \right\|_{\,1}^{\,2}\; \;\forall\; g \,\in\, K\label{eq1.1}
\end{align}
Multiplying (\ref{eq1}) and (\ref{eq1.1}), and using (\ref{eq1.001}) and (\ref{eq1.0001}), we get
\begin{align*}
&A\,C\left\|\,f \,\otimes\, g\,\right\|^{\,2} \,\leq\, \sum\limits_{i,\, j}\,v^{\,2}_{\,i}\,w^{\,2}_{\,j}\,\left<\,\Lambda_{i}\, P_{\,V_{\,i}}\,f,\, \Lambda_{i}^{\,\prime}\,P_{\,V_{\,i}^{\,\prime}}\,f\,\right>_{1}\,\left<\,\Gamma_{j}\,P_{\,W_{\,j}}\,g,\, \Gamma_{j}^{\,\prime}\,P_{\,W_{\,j}^{\,\prime}}\,g\,\right>_{2}\\
& \hspace{2cm}\,\leq\, B\,D\left\|\,f \,\otimes\, g\,\right\|^{\,2}\\
&\Rightarrow\, A\, \left\|\,f \,\otimes\, g\,\right\|^{\,2}\leq\,\sum\limits_{i,\, j}\,v^{\,2}_{\,i}\,w^{\,2}_{\,j}\,\left<\,\Lambda_{i}\, P_{\,V_{\,i}}\,f \,\otimes\, \Gamma_{j}\,P_{\,W_{\,j}}\,g\,,\, \Lambda_{i}^{\,\prime}\,P_{\,V_{\,i}^{\,\prime}}\,f \,\otimes\, \Gamma_{j}^{\,\prime}\,P_{\,W_{\,j}^{\,\prime}}\,g\,\right>\\
&\hspace{2cm}\leq\,B\, \left\|\,f \,\otimes\, g\,\right\|^{\,2}.\\
&\Rightarrow\, A\,C\, \left\|\,f \,\otimes\, g\,\right\|^{\,2}\\
&\leq\,\sum\limits_{i,\, j}\,v^{\,2}_{\,i}\,w^{\,2}_{\,j}\,\left<\left(\,\Lambda_{i} \,\otimes\, \Gamma_{j}\,\right)\,P_{\,V_{\,i} \,\otimes\, W_{\,j}}\,(\,f \,\otimes\, g\,),\, \left(\,\Lambda_{i}^{\,\prime} \,\otimes\, \Gamma_{j}^{\,\prime}\,\right)\,P_{\,V_{\,i}^{\,\prime} \,\otimes\, W_{\,j}^{\,\prime}}\,(\,f \,\otimes\, g\,)\,\right>\\ 
&\hspace{2cm}\leq\,B\,D\,\left\|\,f \,\otimes\, g\,\right\|^{\,2}\; \;\forall\; f \,\otimes\, g \,\in\, H \,\otimes\, K.
\end{align*}
Hence, the family \,$\left(\,\Lambda \,\otimes\, \Gamma,\, \Lambda^{\,\prime} \,\otimes\, \Gamma^{\,\prime}\,\right)$\, is a bi $g$-fusion frame for \,$H \,\otimes\, K$\, with bounds \,$A\,C$\, and \,$B\,D$.\,This completes the proof. 
\end{proof}

\begin{remark}
Let \,$\left(\,\Lambda \,\otimes\, \Gamma,\, \Lambda^{\,\prime} \,\otimes\, \Gamma^{\,\prime}\,\right)$\, be a bi $g$-fusion frame for \,$H \,\otimes\, K$.\,According to the definition (\ref{defn1}), the frame operator \,$S \,:\, H \,\otimes\, K \,\to\, H \,\otimes\, K$\, is described by
\begin{align*}
&S\,(\,f \,\otimes\, g\,) \\
&\,=\sum\limits_{i,\, j}\,v^{\,2}_{\,i}\,w^{\,2}_{\,j}\,P_{\,V_{i}^{\,\prime} \,\otimes\, W_{j}^{\,\prime}}\left(\,\Lambda_{i}^{\,\prime} \,\otimes\, \Gamma_{j}^{\,\prime}\,\right)^{\,\ast}\left(\,\Lambda_{i} \,\otimes\, \Gamma_{j}\,\right)P_{\,V_{i} \,\otimes\, W_{j}}\,(\,f \,\otimes\, g\,)
\end{align*}
for all \,$f \,\otimes\, g \,\in\, H \,\otimes\, K$.  
\end{remark}

\begin{theorem}
If \,$S_{\Lambda,\, \Gamma}, S_{\Lambda^{\,\prime},\, \Gamma^{\,\prime}}$\, and \,$S$\, are the corresponding frame operators for the bi-$g$-fusion frames \,$\left(\,\Lambda,\, \Lambda^{\,\prime}\right)$,\,$\left(\,\Gamma,\,\Gamma^{\,\prime}\,\right)$\, and \,$\left(\,\Lambda \,\otimes\, \Gamma,\, \Lambda^{\,\prime} \,\otimes\, \Gamma^{\,\prime}\,\right)$, respectively, then \,$S \,=\, S_{\Lambda,\, \Gamma} \,\otimes\, S_{\Lambda^{\,\prime},\, \Gamma^{\,\prime}}$\, and \,$S^{\,-\, 1} \,=\, S^{\,-\, 1}_{\Lambda,\, \Gamma} \,\otimes\, S^{\,-\, 1}_{\Lambda^{\,\prime},\, \Gamma^{\,\prime}}$.
\end{theorem}

\begin{proof}
For each \,$f \,\otimes\, g \,\in\, H \,\otimes\, K$, we have
\begin{align*}
&S\,(\,f \,\otimes\, g\,) \\
&\,=\sum\limits_{i,\, j}\,v^{\,2}_{\,i}\,w^{\,2}_{\,j}\,P_{\,V_{i}^{\,\prime} \,\otimes\, W_{j}^{\,\prime}}\left(\,\Lambda_{i}^{\,\prime} \,\otimes\, \Gamma_{j}^{\,\prime}\,\right)^{\,\ast}\left(\,\Lambda_{i} \,\otimes\, \Gamma_{j}\,\right)P_{\,V_{i} \,\otimes\, W_{j}}\,(\,f \,\otimes\, g\,)\\
&=\sum\limits_{i,\, j}v^{\,2}_{\,i}w^{\,2}_{\,j}\left(\,P_{\,V_{i}^{\,\prime}} \,\otimes\, P_{\,W_{j}^{\,\prime}}\,\right)\left(\,\left(\,\Lambda_{i}^{\,\prime}\,\right)^{\,\ast} \,\otimes\, \left(\,\Gamma_{j}^{\,\prime}\,\right)^{\,\ast}\,\right)\left(\,\Lambda_{i} \,\otimes\, \Gamma_{j}\,\right)\left(\,P_{\,V_{i}}\,f \,\otimes\, P_{\,W_{j}}\,g\,\right)\\
&\,=\, \sum\limits_{i,\, j}\,v^{\,2}_{\,i}\,w^{\,2}_{\,j}\,\left(\,P_{\,V_{i}^{\,\prime}}\,\left(\,\Lambda_{i}^{\,\prime}\,\right)^{\,\ast}\,\Lambda_{i}\,P_{\,V_{i}}\,f \,\otimes\, P_{\,W_{j}^{\,\prime}}\,\left(\,\Gamma_{j}^{\,\prime}\,\right)^{\,\ast}\,\Gamma_{j}\,P_{\,W_{j}}\,g\,\right)\\
&\,=\, \left(\,\sum\limits_{\,i \,\in\, I}\,v^{\,2}_{\,i}\,P_{\,V_{i}^{\,\prime}}\,\left(\,\Lambda_{i}^{\,\prime}\,\right)^{\,\ast}\,\Lambda_{i}\,P_{\,V_{i}}\,f\,\right) \,\otimes\, \left(\,\sum\limits_{\,j \,\in\, J}\,w^{\,2}_{j}\,P_{\,W_{j}}\,\left(\,\Gamma_{j}^{\,\prime}\,\right)^{\,\ast}\,\Gamma_{j}\,P_{\,W_{j}}\,g\,)\,\right)\\
&=\,S_{\Lambda,\, \Gamma}\,f \,\otimes\, S_{\Lambda^{\,\prime},\, \Gamma^{\,\prime}}\,g \,=\, \left(\,S_{\Lambda,\, \Gamma} \,\otimes\, S_{\Lambda^{\,\prime},\, \Gamma^{\,\prime}}\,\right)\,(\,f \,\otimes\, g\,) 
\end{align*} 
This shows that \,$S \,=\,S_{\Lambda,\, \Gamma} \,\otimes\, S_{\Lambda^{\,\prime},\, \Gamma^{\,\prime}}$.\,Since \,$S_{\Lambda,\, \Gamma}$\, and \,$S_{\Lambda^{\,\prime},\, \Gamma^{\,\prime}}$\, are invertible, by \,$(\,iv\,)$\, of Theorem \ref{th1.1}, it follows that \,$S^{\,-\, 1} \,=\, \left(\,S_{\Lambda,\, \Gamma} \,\otimes\, S_{\Lambda^{\,\prime},\, \Gamma^{\,\prime}}\,\right)^{\,-\, 1} \,=\, S^{\,-\, 1}_{\Lambda,\, \Gamma} \,\otimes\, S^{\,-\, 1}_{\Lambda^{\,\prime},\, \Gamma^{\,\prime}}$.\,This completes the proof.     
\end{proof}

\begin{proposition}
If \,$S_{\Lambda,\, \Gamma}, S_{\Lambda^{\,\prime},\, \Gamma^{\,\prime}}$\, and \,$S$\, are the corresponding frame operators for the bi-$g$-fusion frames \,$\left(\,\Lambda,\, \Lambda^{\,\prime}\right)$,\,$\left(\,\Gamma,\,\Gamma^{\,\prime}\,\right)$\, and \,$\left(\,\Lambda \,\otimes\, \Gamma,\, \Lambda^{\,\prime} \,\otimes\, \Gamma^{\,\prime}\,\right)$, respectively.\,Then
\[A\,C\,I_{H\, \otimes\, K} \,\leq\, S \,\leq\,  B\,D\,I_{H\, \otimes\, K},\]
where \,$(\,A,\,B\,)$\, and \,$(\,C,\, D\,)$\, are frame bounds of \,$\Lambda_{T\,U}$\, and \,$\Gamma_{T_{1}\,U_{1}}$, respectively and \,$I_{H\, \otimes\, K}$\, is the identity operator on \,$H \,\otimes\, K$. 
\end{proposition}

\begin{proof}
Suppose \,$\left(\,\Lambda,\, \Lambda^{\,\prime}\right)$\, and \,$\left(\,\Gamma,\,\Gamma^{\,\prime}\,\right)$\, are bi-$g$-fusion frames for \,$H$\, and \,$K$\, with bounds \,$(\,A,\, B\,)$\, and \,$(\,C,\, D\,)$, respectively.\,Then 
\[A\,I_{H} \,\leq\, S_{\Lambda,\, \Gamma} \,\leq\, B\,I_{H},\, \,C\,I_{K} \,\leq\, S_{\Lambda^{\,\prime},\, \Gamma^{\,\prime}} \,\leq\, D\,I_{K}.\]
Taking tensor product of the above two inequalities we get
\begin{align*}
&A\,C\,\left(\,I_{H} \,\otimes\, I_{K}\,\right) \,\leq\, S_{\Lambda,\, \Gamma} \,\otimes\, S_{\Lambda^{\,\prime},\, \Gamma^{\,\prime}} \,\leq\, B\,D\,\left(\,I_{H} \,\otimes\, I_{K}\,\right)\\
&\Rightarrow\,A\,C\,I_{H\, \otimes\, K} \,\leq\, S \,\leq\,  B\,D\,I_{H\, \otimes\, K}.
\end{align*}
\end{proof}

\begin{theorem}\label{th2}
Suppose \,$U_{1} \,\in\, \mathcal{B}\,(\,H\,)$\, and \,$U_{2} \,\in\, \mathcal{B}\,(\,K\,)$.\,Let \,$\left(\,\Lambda,\, \Lambda^{\,\prime}\right)$\, and \,$\left(\,\Gamma,\,\Gamma^{\,\prime}\,\right)$\, be bi-$g$-fusion frames for \,$H$\, and \,$K$\, with bounds \,$A,\, B$\, and \,$C,\, D$, respectively having their corresponding frame operators \,$S_{\Lambda,\, \Gamma}, S_{\Lambda^{\,\prime},\, \Gamma^{\,\prime}}$, respectively.\,Then the pair \,$\left\{\,\left(\,\left(\,U_{1} \,\otimes\, U_{2}\,\right)\left(\,V_{i} \,\otimes\, W_{j}\,\right),\, \left(\,\Lambda_{i} \,\otimes\, \Gamma_{j}\,\right)P_{\,V_{i} \,\otimes\, W_{j}}\,\left(\,U_{1} \,\otimes\, U_{2}\,\right)^{\,\ast},\, v_{i}\,w_{j}\,\right)\,\right\}_{\,i,\,j}$\, and \,$\left\{\,\left(\,\left(\,U_{1} \,\otimes\, U_{2}\,\right)\left(\,V_{i}^{\,\prime} \,\otimes\, W_{j}^{\,\prime}\,\right),\, \left(\,\Lambda_{i}^{\,\prime} \,\otimes\, \Gamma_{j}^{\,\prime}\,\right)P_{\,V_{i}^{\,\prime} \,\otimes\, W_{j}^{\,\prime}}\,\left(\,U_{1} \,\otimes\, U_{2}\,\right)^{\,\ast},\, v_{i}\,w_{j}\,\right)\,\right\}_{\,i,\,j}$\, bi-$g$-fusion frame for \,$H \,\otimes\, K$.   
\end{theorem}

\begin{proof}
Take \,$\Theta_{i\,j} \,=\, \left(\,\Lambda_{i} \,\otimes\, \Gamma_{j}\,\right)\,P_{\,V_{i} \,\otimes\, W_{j}}\,\left(\,U_{1} \,\otimes\, U_{2}\,\right)^{\,\ast}\,P_{\left(\,U_{1} \,\otimes\, U_{2}\,\right)\left(\,V_{i} \,\otimes\, W_{j}\,\right)}$\, and \,$\Theta_{i\,j}^{\,\prime} \,=\, \left(\,\Lambda_{i}^{\,\prime} \,\otimes\, \Gamma_{j}^{\,\prime}\,\right)\,P_{\,V_{i}^{\,\prime} \,\otimes\, W_{j}^{\,\prime}}\,\left(\,U_{1} \,\otimes\, U_{2}\,\right)^{\,\ast}\,P_{\left(\,U_{1} \,\otimes\, U_{2}\,\right)\left(\,V_{i}^{\,\prime} \,\otimes\, W_{j}^{\,\prime}\,\right)}$.\,Then using Theorem \ref{th0.001} and Theorem \ref{th1.1}, we get
\begin{align}
\Theta_{i\,j} &\,=\, \left(\,\Lambda_{i} \,\otimes\, \Gamma_{j}\,\right)\,\left(\,P_{\,V_{i}} \,\otimes\, P_{W_{j}}\,\right)\,\left(\,U_{1}^{\,\ast} \,\otimes\, U_{2}^{\,\ast}\,\right)\,\left(\,P_{\,U_{1}\,V_{i}} \,\otimes\, P_{U_{2}\,W_{j}}\,\right)\nonumber\\
&=\,\Lambda_{i}\,P_{\,V_{i}}\,U_{1}^{\,\ast}\,P_{\,U_{1}\,V_{i}} \,\otimes\, \Gamma_{j}\,P_{W_{j}}\,U_{2}^{\,\ast}\,P_{U_{2}\,W_{j}} \,=\, \Lambda_{i}\,P_{\,V_{i}}\,U_{1}^{\,\ast} \,\otimes\, \Gamma_{j}\,P_{W_{j}}\,U_{2}^{\,\ast}\label{eqq1}  
\end{align}
\begin{align}
\Theta_{i\,j}^{\,\prime} =\,\Lambda_{i}^{\,\prime}\,P_{\,V_{i}^{\,\prime}}\,U_{1}^{\,\ast}\,P_{\,U_{1}\,V_{i}^{\,\prime}} \,\otimes\, \Gamma_{j}^{\,\prime}\,P_{W_{j}^{\,\prime}}\,U_{2}^{\,\ast}\,P_{U_{2}\,W_{j}^{\,\prime}} \,=\, \Lambda_{i}^{\,\prime}\,P_{\,V_{i}^{\,\prime}}\,U_{1}^{\,\ast} \,\otimes\, \Gamma_{j}^{\,\prime}\,P_{W_{j}^{\,\prime}}\,U_{2}^{\,\ast}\label{eqqn1}  
\end{align}
Now, using (\ref{eqq1}) and (\ref{eqqn1}), for each \,$f \,\otimes\, g \,\in\, H \,\otimes\, K$, we have
\begin{align*}
&\sum\limits_{i,\, j}\,v^{\,2}_{\,i}\,w^{\,2}_{\,j}\,\left<\,\Theta_{i\,j}\,(\,f \,\otimes\, g\,),\, \Theta_{i\,j}^{\,\prime}\,(\,f \,\otimes\, g\,)\,\right>\\
&=\,\sum\limits_{i,\, j}\,v^{\,2}_{\,i}\,w^{\,2}_{\,j}\,\left<\,\Lambda_{i}\,P_{\,V_{i}}\,U_{1}^{\,\ast}\,f \,\otimes\, \Gamma_{j}\,P_{W_{j}}\,U_{2}^{\,\ast}\,g\,,\, \Lambda_{i}^{\,\prime}\,P_{\,V_{i}^{\,\prime}}\,U_{1}^{\,\ast}\,f \,\otimes\, \Gamma_{j}^{\,\prime}\,P_{W_{j}^{\,\prime}}\,U_{2}^{\,\ast}\,g\,\right>\\
&=\,\sum\limits_{i,\, j}\,v^{\,2}_{\,i}\,w^{\,2}_{\,j}\,\left<\,\Lambda_{i}\,P_{\,V_{i}}\,U_{1}^{\,\ast}\,f,\, \Lambda_{i}^{\,\prime}\,P_{\,V_{i}^{\,\prime}}\,U_{1}^{\,\ast}\,f\,\right>_{1}\,\left<\, \Gamma_{j}\,P_{W_{j}}\,U_{2}^{\,\ast}\,g,\, \Gamma_{j}^{\,\prime}\,P_{W_{j}^{\,\prime}}\,U_{2}^{\,\ast}\,g\,\right>_{2}\\
&=\,\sum\limits_{i \,\in\, I}\,v^{\,2}_{\,i}\,\left<\,\Lambda_{i}\,P_{\,V_{i}}\,U_{1}^{\,\ast}\,f,\, \Lambda_{i}^{\,\prime}\,P_{\,V_{i}^{\,\prime}}\,U_{1}^{\,\ast}\,f\,\right>_{1}\sum\limits_{j \,\in\, J}\,w^{\,2}_{\,j}\,\left<\, \Gamma_{j}\,P_{W_{j}}\,U_{2}^{\,\ast}\,g,\, \Gamma_{j}^{\,\prime}\,P_{W_{j}^{\,\prime}}\,U_{2}^{\,\ast}\,g\,\right>_{2}\\
&\leq\, B\,\left\|\,U_{1}^{\,\ast}\,f\,\right\|^{\,2}_{1}\,D\,\left\|\,U_{2}^{\,\ast}\,g\,\right\|^{\,2}_{2}\leq\, B\,D\,\left\|\,U_{1}\,\right\|^{\,2}\,\|\,f\,\|^{\,2}_{1}\,\left\|\,U_{2}\,\right\|^{\,2}\,\|\,g\,\|_{2}\\
&=\,B\,D\,\left\|\,U_{1} \,\otimes\, U_{2}\,\right\|^{\,2}\,\|\,f \,\otimes\, g\,\|^{\,2}. 
\end{align*}
On the other hand, for each \,$f \,\otimes\, g \,\in\, H \,\otimes\, K$, we have
\begin{align*}
&\sum\limits_{i,\, j}\,v^{\,2}_{\,i}\,w^{\,2}_{\,j}\,\left<\,\Theta_{i\,j}\,(\,f \,\otimes\, g\,),\, \Theta_{i\,j}\,(\,f \,\otimes\, g\,)\,\right>\\
&=\,\sum\limits_{i \,\in\, I}\,v^{\,2}_{\,i}\,\left<\,\Lambda_{i}\,P_{\,V_{i}}\,U_{1}^{\,\ast}\,f,\, \Lambda_{i}^{\,\prime}\,P_{\,V_{i}^{\,\prime}}\,U_{1}^{\,\ast}\,f\,\right>_{1}\sum\limits_{j \,\in\, J}\,w^{\,2}_{\,j}\,\left<\, \Gamma_{j}\,P_{W_{j}}\,U_{2}^{\,\ast}\,g,\, \Gamma_{j}^{\,\prime}\,P_{W_{j}^{\,\prime}}\,U_{2}^{\,\ast}\,g\,\right>_{2}\\
&\geq\, A\,\left\|\,U_{1}^{\,\ast}\,f\,\right\|^{\,2}_{1}\,C\,\left\|\,U_{2}^{\,\ast}\,g\,\right\|^{\,2}_{2}\geq\, A\,C\,\left\|\,U_{1}^{\,-\,1} \,\right\|^{\,-\,2}\,\left\|\,U_{2}^{\,-\,1} \,\right\|^{\,-\,2}\,\|\,f\,\|^{\,2}_{1}\,\|\,g\,\|_{2}\\
&=\,A\,C\,\left\|\,\left(\,U_{1} \,\otimes\, U_{2}\,\right)^{\,-\, 1}\,\right\|^{\,-\,2}\,\|\,f \,\otimes\, g\,\|^{\,2}. 
\end{align*}
This completes the proof.
\end{proof}

\end{document}